\documentclass{article}
\usepackage{amsfonts,amssymb,amsmath,amscd,amsthm}
\usepackage{graphicx}
\usepackage{epstopdf}
\usepackage{mathrsfs}
\usepackage{enumitem}

\date{}
\newcommand\be{\begin{array}}
\newcommand\en{\end{array}}
\newcommand\di{\displaystyle}

\newcommand\ga{\gamma}

\newcommand\la{\lambda}
\newcommand\Om{\Omega}
\newcommand\wi{\widetilde}

\newcommand\gl{\geqslant}
\newcommand\ls{\leqslant}

\newcommand\ri{\rightarrow}

\newcommand\bk{\backslash}

\newcommand\mb{\mathbb}
\newcommand\ov{\overline}

\renewcommand{\emptyset}{\varnothing}

\begin{document}
\title{\bf      Solutions  for Strongly Monotone Operator Equations in Riesz Spaces
\thanks{This paper is supported   by the National Natural Science Foundation of China (Grant No. 11871250).}
}
\author
{Xian  Xu$^1$,  Baoxia Qin$^2$} \maketitle

\vspace{2mm} \noindent \begin{center}\ $^1$Department of Mathematics,
Jiangsu Normal University, Xuzhou,
\\Jiangsu,
221116,  P. R. China

$^2$School of Mathematics,
Qilu Normal  University, Jinan,
\\Shandong,
250013,  P. R. China\\

\end{center}
\begin{center}
\begin{minipage}{5in}
{\small {\bf Abstract}\quad This paper is devoted to the study of solutions for a class of operator equations governed by strongly monotone operators on real Riesz spaces continuously embedded into Banach lattices.
By exploiting the intrinsic lattice structure of Banach lattices, we establish refined growth assumptions on nonlinear terms to guarantee that suitable neighbourhoods of positive and negative cones are invariant under the descending flow. Combining the descending flow invariant set technique with the theory of strongly monotone operators, we derive abstract existence theorems: the operator equation possesses at least one positive solution, one negative solution and one sign-changing solution. These abstract results are further applied to \((p,q)\)-Laplacian boundary value problems, yielding corresponding multiplicity conclusions on positive, negative and sign-changing solutions.


}
\end{minipage}
\end{center}

\section{I. Introduction}

Riesz spaces have been systematically investigated in linear functional analysis, with comprehensive fundamental theories well documented in standard monographs, see [1]. Over the past two decades, researchers have extended lattice methods to nonlinear functional analysis on Riesz spaces. For instance, Sun and Liu [2,3] combined lattice order structures with topological degree theory to obtain fixed-point existence results for nonlinear operators. Sun and Xu [4,5] integrated lattice theory with bifurcation analysis to characterize the global structure of solution branches for parameter-dependent nonlinear operator equations.

Critical point theory on ordered vector lattices has also been developed by many authors. Perera and Schechter [6] employed variational methods on Riesz spaces to treat Fu\v{c}\'{i}k spectra and jumping nonlinearities. Kajikiya [7] incorporated lattice decomposition into deformation lemmas and established a mountain pass principle in ordered Banach spaces, with applications to semilinear elliptic boundary value problems.

The primary purpose of the present paper is to further develop critical point theory for Riesz spaces. More precisely, we establish existence of sign-changing, positive and negative critical points for equations driven by strongly monotone operators via the descending flow invariant set method.

The invariant set method for descending flow was originally proposed by Sun Jingxian in [9], and has become a standard variational tool for elliptic boundary value problems. Its core advantage is that the search for distinct critical points reduces to constructing appropriate flow-invariant closed sets with nonempty interior. Numerous works have adopted this strategy to study multiple solutions of elliptic PDEs; we refer to [8,9,10,15].
and proved multiplicity of solutions by constructing descending flow invariant subsets in $W_0^{1,p}(\Omega)$.

A well-known technical obstacle arises when applying this method: standard order intervals in Sobolev spaces generically have empty interior, so the pair of conditions (Palais--Smale compactness + closed convex invariant sets with interior) rarely hold simultaneously. To overcome this difficulty, Bartsch--Liu [8] and Liu--Sun [10] adopted a two-space framework, working simultaneously in $W_0^{1,p}(\Omega)$ and $C_0^1(\Omega)$. An alternative approach was proposed in [11], replacing closed convex cones with their open neighbourhoods, which automatically possess nonempty interior.

This naturally motivates the central question of our work: can we fully incorporate Banach lattice order structures into the descending flow invariant set argument to obtain sign-changing critical points directly on the underlying Banach space?

We address this question by studying operator equations on reflexive Banach spaces. Let $X$ be a real reflexive Banach space with dual $X^*$. We treat equations of the form
\[
T(u) = \mathbf{f}(u),\tag{1.1}
\]
where $T:X\to X^*$ is a strongly monotone operator and $\mathbf{f}:X\to X^*$ is a nonlinear operator.

Inspired by the two-space variational scheme in [11], we combine monotone operator theory on reflexive Banach spaces with descending flow invariant sets to prove existence of positive, negative and sign-changing solutions for (1.2). Strongly monotone operators arise naturally in many elliptic boundary value problems, including those governed by $p$-Laplacian and $(p,q)$-Laplacian differential operators; see [12,13,14,16,17,18]. Consequently, our abstract variational theorems apply directly to such nonlinear elliptic PDEs.

\section{Main Results}

\ \ \ \ \ \rm  First let us recall some concepts on Riesz space. Let \(X, E\) be   real reflexive Banach spaces with the norm $\|\cdot\|$ and $\|\cdot\|_1$, respectively. And let \(X^{*}\) and \(E^*\) be the conjugate spaces of \(X\) and $E$, respectively. Let \(\langle\cdot,\cdot\rangle\) denote the conjugate pair of \(X^{*}\) and \(X\), \((\cdot,\cdot)\) denote the conjugate pair of \(E^{*}\) and \(E\), respectively. Assume that $X\hookrightarrow\hookrightarrow E$. Then $\|x\|_1\ls c_0\|x\|$ for each $x\in X$ and some $c_0>0$.  Let $ P_1$ be a cone of $E$ which induce an ordering $\ls$, that is $x\ls y$ if and only if $y-x\in P_1$.  Obviously, $P=P_1\cap X$ is also a cone in $X$ which induce a ordering in $X$, denoted also by $\ls$. Let $P_1^*$ be the dual cone of $P_1$, i.e.,
$$P_1^*=\{g\in E^*: (g, u)\gl 0\ \mbox{for all}\ u\in P_1\}.$$

 The pair $(E, \ls)$ is assumed to be a  Riesz space. This implies that $\ls$ is compatible with the algebraic structure of $E$, i.e.,

(1)\ For any $x, y, z\in E$ with $x\ls y$, $x+z\ls y+z$,

(2)\ For any $x\in E$ with $x\gl 0$,  and $r\gl 0$,  $rx\gl 0$,\\
and every set $\{x, y\}\subset E$ has a  supremum  $\sup\{x,y\}=x\vee y$.

Let $x\in E$. The elements
$$x^+=x\vee 0,\  x^-=(-x)\vee 0,\  |x|=x^++x^-$$
are called the positive part, negative part, and modulus of  $x$, respectively. Obviously, $0\ls x^+\ls |x|$ and $0\ls x^-\ls |x|$. Moreover,
 $x\ls y$ if and only if $y^+\gl x^+$ and $x^-\gl y^-$. In what follows we write $x_+=x^+$ and $x_-=-x^-$. Then, we have $x= x_++x_-$. We additionally  assume that $(E, \ls)$ is  a Banach lattice, and so the norm $\|\cdot\|_1$ is  a Riesz norm, i.e., $|x|\ls |y|$ implies  $\|x\|_1\ls \|y\|_1$.
\medskip

{\bf Definition 2.1([1])}\  Let $E$ be a Riesz space, $x, y\in E$. Then $x$ and $y$ are said to be  disjoint  if $|x|\wedge|y|=0$, denoted by $x\bot y$.
\medskip

{\bf Definition 2.2([1])}\  Let $E$ be a Riesz space,  $F$ be a real linear space and let the operator $T: E\to F$. Then
  $T$ is said to be  orthogonally additive if $T(x+y)=T(x)+T(y)$ for all $x, y\in E$ with $x\bot y$.
\medskip

 In this paper, we make the following  assumptions:

(H\(_{1})\)\ \(J:X\to \mathbb R\) is of \( C^1\), and \(J^{\prime}(u)\) has the following form
\[\langle J^{\prime}(u),v\rangle = \langle T(u),v\rangle-(\mathbf{f}(u),v), \forall v\in X,\]
where \(\mathbf{f}:E\to E^*\) is continuous and \(T:X\to X^{*}\) is continuous, \(\langle T(u), u_{\pm}\rangle=\langle T(u_\pm), u_{\pm}\rangle\)  for any \(u\in X\), and
\[\langle T(u) - T(v),u - v\rangle\gl c_{1}\|u - v\|^{p},\ \forall u,v\in X\eqno(2.1)\]
for some \(c_{1}>0\) and \(p\gl 2\).

 (H\(_2)\)\ \(\mathbf{f}:E\to E^{*}\) is orthogonally additive,  $\mathbf{f}(\pm P)\subset \pm P^*_1$ and
 $$( \mathbf{f}(u), w)\ls g(\|u\|_1)\|w\|_1$$ for each $u,w\in X$, where  $g: \mb R^+\ri \mb R^+$ is increasing and $\lim\limits_{s\ri 0^+}\frac{g(s)}{s^{p-1}}=0$.

(H\(_3)\)\ For each $u\in X\bk\{0\}$, $\lim\limits_{t\ri+\infty}J(tu)=-\infty$.

\medskip

The following Definition 2.3, Lemma 2.1 and Corollary 2.1  can be found in [14].
\medskip

{\bf Definition 2.3}\
Let \(D\subset X\) and the mapping \(T:D\to X^{*}\). Let \(x_{0}\in D\). If \(h\in X\), \(t_{n}>0\), \(x_{0}+t_{n}h\in D\), \(t_{n}\to 0\Rightarrow T(x_{0}+t_{n}h)\rightharpoondown Tx_{0}\), then
\(T\) is said to be hemicontinuous at \(x_{0}\). If \(T\) is hemicontinuous at every point in \(D\), then \(T\) is said to be hemicontinuous on \(D\).
\medskip

{\bf Lemma 2.1}\ ( Strongly Monotone Mapping)\ Let \(X\) be reflexive and the mapping \(T: X \to X^{*}\) be hemicontinuous and satisfy the strongly monotone condition
\[
\langle T(x) - T(y), x - y\rangle \gl \alpha(\|x - y\|)\|x - y\|, \forall x, y \in X,\eqno(2.2)
\]
where \[\alpha(0) = 0, \ \alpha(t) > 0\ (\forall t > 0);\  \lim\limits_{t \to +\infty} \alpha(t)=+\infty,\]
then \(T\) is surjective and one-to-one (that is, for any \(\wi g\in X^{*}\), the equation \(T(x) =\wi g\) has a unique solution in \(X\)).
\medskip

{\bf Corollary 2.1}\  Suppose that all the conditions of Lemma 2.1 are satisfied and \(T: X \to X^{*}\) is continuous, and the function \(\alpha(t)\) is continuous on \(0 < t < +\infty\). Then \(T\) is a homeomorphism between \(X\) and \(X^{*}\).
\medskip

From Corollary 2.1, if (H\(_{1})\) holds, then there exists a mapping \(A: X \to X\) such that for any \(u \in X\), \(v = A(u)\), satisfying
\[
\langle T(v), w\rangle = (\mathbf{f}(u), w),\  \forall w \in X.
\]
And since \(\mathbf{f}: E \to E^*\subset X^*\) is continuous, we know that  \(A: X \to X\) is continuous.
\medskip

 Obviously, from (H$_1$), \(u\) being a solution of (1.2) is equivalent to \(u\) being a critical point of \(J\). Let \(K = \{u\in X: J'(u) = 0\}\). For any \(u\in K\setminus\{0\}\), we call \(u\) a positive solution of (1.2) when \(u\in P\), a negative solution when \(u\in -P\), and a sign-changing solution when \(u\notin P\cup(-P)\).
\medskip

We introduce a compactness condition replacing the standard Palais-Smale  (P.S.) condition.
\begin{itemize}
\item[(H\(_4)\)]\ Let $\{u_n\}\subset X$ satisfy $\{u_n\}\subset J^{-1}([a,b])$ for some bounded interval $[a,b]\subset\mb R$ and $\|u_n-A(u_n)\|\to 0$. Then $\{u_n\}$ is bounded in $X$ and $J'(u_n)\to 0$ in $X^*$ as $n\to\infty$.
\end{itemize}
Now we have the following main result:

{\bf Theorem 2.1}\ Suppose  (H\(_{1})\sim\)(H\(_{4})\) hold.
Then the operator equation \((1.1)\) admits at least one positive, one negative and one sign-changing solution.
\medskip

{\bf Remark 2.1}\  The  growth bound (2.1) can be replaced by the more general strongly monotone structure (2.2) without altering the proof logic.
\medskip

We now develop auxiliary lemmas required for Theorem 2.1. Throughout this section we always assume $(\mathrm{H}_1)$--$(\mathrm{H}_4)$ are satisfied.\medskip

{\bf Lemma 2.2}\ For any given \(\mu> 0\), define the open neighbourhoods of the cones \[D^{\pm}(\mu)=\{x \in X: \mbox{dist}\ (x, \pm P)<\mu\}.\] There exists \(\mu_{0} > 0\) such that for any \(0 < \mu < \mu_{0}\), we have \(A(D^{\pm}(\mu)) \subset D^{\pm}(\frac{1}{2}\mu)\).

{\bf Proof.}\
Fix $u\in X$, set $v=A(u)$. Since $\mathbf{f}(P)\subset P^*_1$ and $v_-\in (-P)\subset (-P_1)$, if $u\in P$ we have
\[c_1\|v_-\|^p\ls\langle T(v_-), v_-\rangle=\langle T(v), v_-\rangle = ( \mathbf{f}(u), v_-)\ls 0,\]
and so $v_-=0$, $v=v_+\in P$. This yields that $A(P)\subset P$.

By (H$_1)$ and (H$_2)$ we have for any $u\in D^+(\mu)$,
\begin{align*}
c_1\|v_-\|^p&\ls\langle T(v_-), v_-\rangle=\langle T(v), v_-\rangle\\
  &= ( \mathbf{f}(u), v_-)=( \mathbf{f}(u_+)+\mathbf{f}(u_-), v_-)\\
  &\ls ( \mathbf{f}(u_-), v_-)\ls  g(\|u_-\|_1)\|v_-\|_1\\
  &\ls  c_0g(\|u_-\|_1)\|v_-\|,\\
\end{align*}
and so $c_1\|v_-\|^{p-1}\ls c_0g(\|u_-\|_1)$.  Since $v=v_++v_-$ and $v_+\in P$, we have
\[c_1\big(\mbox{dist}\ (v, P)\big)^{p-1}\ls c_1\|v_-\|^{p-1}\ls c_0g(\|u_-\|_1).\]
For each $w\in P$, we have $u-w\ls u$, and so $$|u_-|=u^-\ls (u-w)^-\ls |u-w|.$$  Since  $\|\cdot\|_1$ is a Riesz norm and $g: \mb R^+\ri \mb R^+$ is increasing, we have
\[c_1\big(\mbox{dist}\ (v, P)\big)^{p-1}\ls c_0 g(\|u_-\|_1)\ls c_0g(\|u-w\|_1)\ls c_0g(c_0\|u-w\|).\]
Consequently, we have $$c_1\big(\mbox{dist}\ (v, P)\big)^{p-1}\ls c_0g(c_0\mbox{dist}\ (u, P)).$$ Since $\lim\limits_{s\ri 0^+}\frac{g(s)}{s^{p-1}}=0$, there exists \(\mu_{0}^+ > 0\) such that for any \(0 < \mu < \mu_{0}^+\), we have  \(A(D^{+}(\mu)) \subset D^{+}(\frac{1}{2}\mu)\).

Similarly, there exists \(\mu_{0}^- > 0\) such that for any \(0 < \mu < \mu_{0}^-\), we have \(A(D^{-}(\mu)) \subset D^{-}(\frac{1}{2}\mu)\). Let $\mu_0=\min\{\mu^+_0,\mu^-_0\}$. Then the conclusion holds. The proof is complete.
\medskip

Let \(\mu\in(0,\mu_{0})\) be given. For brevity, write \(D_{1}=D^{+}(\mu), D_{2}=D^{-}(\mu), D_{3}=D_{1}\cap D_{2}\).

{\bf Lemma 2.3}\ Denote \(\widetilde{X}:=X\backslash K\). There exists a locally Lipschitz mapping \(B:\widetilde{X}\to X\) satisfying

(1)\   \(B(D_1)\subset D^{+}(\frac{1}{2}\mu)\),  \(B(D_2)\subset D^{-}(\frac{1}{2}\mu)\);

(2)\ \(\frac{1}{2}\|u - A(u)\|\ls\|u - B(u)\|\ls 2\|u - A(u)\|, \forall u\in\widetilde{X}\);

 (3)\ \(\langle J^{\prime}(u),u - B(u)\rangle\gl\frac{c_{1}}{2}\|u - A(u)\|^{p}, \forall u\in\widetilde{X}\).

{\bf Proof.}\
Write \(\overline{D}_{i}(i = 1,2,3)\) to represent the closure of \(D_{i}\) in \(X\).
For any given \(u_{0}\in\widetilde{X}\backslash(\overline{D}_{1}\cup\overline{D}_{2})\), choose a neighborhood \(U_1(u_{0}))\) of $u_0$ such that
\(U_1(u_{0}))\subset\widetilde{X}\backslash(\overline{D}_{1}\cup\overline{D}_{2})\). Due to the continuity of \(A\), we can further require that,
for any \(v,w\in U_1(u_0)\),
\[\|A(v) - A(u_{0})\|<\frac{1}{8}\|u_{0} - A(u_{0})\|<\frac{1}{4}\|v - A(v)\|,\eqno(2.3)\]
\[\|A(w) - A(v)\|<\frac{1}{8}\big(\|J'(u_{0})\|\big)^{-1}\|u_{0} - A(u_{0})\|^{p}, \eqno(2.4)\]
\[\|J'(u_0)\|\gl \frac{1}{2}\|J'(v)\|.\eqno(2.5)\]
Let \(\mathcal{A}_{1}=\{U_1(u_{0}):u_{0}\in\widetilde{X}\backslash(\overline{D}_{1}\cup\overline{D}_{2})\}.\)

For any \(u_{0}\in(\widetilde{X}\cap\overline{D}_{1})\backslash\overline{D}_{2}\),  choose a  neighborhood \(U_2(u_{0})\) of $u_0$ such that \(U_2(u_{0})\subset \widetilde{X}\backslash\overline{D}_{2}\),
  and \((2.3)\)$\sim$\((2.5)\) hold for any \(v,w\in U_2(u_{0})\). Let
\(\mathcal{A}_{2}=\{U_{2}(u_{0}):u_{0}\in(\widetilde{X}\cap\overline{D}_{1})\backslash\overline{D}_{2}\}\).

For any \(u_{0}\in(\widetilde{X}\cap\overline{D}_{2})\backslash\overline{D}_{1}\),  choose a  neighborhood \(U_3(u_{0})\) of $u_0$ such that \(U_3(u_{0})\subset \widetilde{X}\backslash\overline{D}_{1}\),
  and \((2.3)\)$\sim$\((2.5)\) hold for any \(v,w\in U_3(u_{0})\). Let
\(\mathcal{A}_{3}=\{U_{3}(u_{0}):u_{0}\in(\widetilde{X}\cap\overline{D}_{2})\backslash\overline{D}_{1}\}\).

For any \(u_{0}\in\wi X\cap\overline{D_{3}}\),  choose a  neighborhood \(U_4(u_{0})\) of $u_0$ such that \(U_4(u_{0})\subset\widetilde{X}\), and \((2.3)\)$\sim$\((2.5)\) hold for any \(v,w\in U_4(u_{0})\). Let
\(\mathcal{A}_{4}=\{U_{4}(u_{0}):u_{0}\in\wi X\cap\overline{D}_{3}\}\).

Finally, let \(\mathcal{A}=\bigcup\limits_{i=1}^4\mathcal{A}_{i}\). Obviously, \(\mathcal{A}\) is an open cover of \(\widetilde{X}\). By using the paracompactness of \(\widetilde{X}\)
we know that there exists a locally finite refinement cover \(\{V_{\lambda}:\lambda\in\Lambda\}\) of \(\mathcal{A}\). Let
\(\mu_{\lambda}\), $\la\in \Lambda$, be a locally finite partition of unity subordinated to \(\{V_{\lambda}:\lambda\in\Lambda\}\).

For any \(\lambda\in\Lambda\), since \(\{V_{\lambda}:\lambda\in\Lambda\}\) is a refinement cover of \(\mathcal{A}\), there must exist \(u_{\lambda}\in\widetilde{X}\) and
\(i\in\{1,2,3,4\}\) such that \(V_{\lambda}\subset U_{i}(u_{\lambda})\). For any \(\lambda\in\Lambda\), we take \(z_\lambda\) as follows:
\begin{itemize}
\item If \(V_{\lambda}\cap\overline{D}_{3}\neq\emptyset\), take any \(z_\lambda\in V_{\lambda}\cap\overline{D}_{3}\);
\item If \(V_{\lambda}\cap(\overline{D}_{1}\backslash\overline{D}_{3})\neq\emptyset\), take any \(z_\lambda\in V_{\lambda}\cap(\overline{D}_{1}\backslash\overline{D}_{3})\);
\item If \(V_{\lambda}\cap(\overline{D}_{2}\backslash\overline{D}_{3})\neq\emptyset\), take any \(z_\lambda\in V_{\lambda}\cap(\overline{D}_{2}\backslash\overline{D}_{3})\);
\item In other cases, take any \(z_\lambda\in V_{\lambda}\).
\end{itemize}

For any \(u\in\widetilde{X}\), let
\[
B(u)=\sum_{\lambda\in\Lambda}\mu_{\lambda}(u)A(z_{\lambda}).\eqno(2.6)
\]
Clearly, \(B:\widetilde{X}\to X\) is a locally Lipschitz continuous mapping, so \(I - B\) is locally Lipschitz continuous. Since $A(D^\pm(\mu))\subset D^\pm(\frac{1}{2}\mu)$, by (2.6) we easily see that $B(D^\pm(\mu))\subset D^\pm(\frac{1}{2}\mu)$.

For any \(u\in\widetilde{X}\), due to the local finiteness of \(\{V_{\lambda}:\lambda\in\Lambda\}\), there are only finitely many sets in \(\{V_{\lambda}:\lambda\in\Lambda\}\), denoted as \(V_{\lambda_{i}}\ (i = 1,2,\cdots,n(u))\), satisfying \(u\in V_{\lambda_{i}}\). For any \(\lambda_{i}\), there exists \(U_{j}(u_{\lambda_{i}})\in\mathcal{A}\ (j\in\{1,2,3,4\})\)
such that \(V_{\lambda_{i}}\subset U_{j}(u_{\lambda_{i}})\). Then, for \(i = 1,2,\cdots,n(u)\), by \((2.3)\)$\sim$\((2.5)\), we know that
\[\be{ll}
\|A(u) - A(z_{\lambda_{i}})\|&\ls\|A(u) - A(u_{\lambda_{i}})\|+\|A(u_{\lambda_{i}}) - A(z_{\lambda_{i}})\|\\
&<\frac{1}{4}\|u_{\lambda_{i}} - A(u_{\lambda_{i}})\|\ls\frac{1}{2}\|u - Au\|.
\en\]
Consequently, by  (2.6) we have for any \(u\in\widetilde{X}\),
\[
\|A(u) - B(u)\|\ls\sum_{i=1}^{n(u)}\mu_{\lambda_i}(u)\|A(u) - A(z_{\lambda_i})\|\ls \frac{1}{2}\|u - Au\|,\] and thus \[\frac{1}{2}\|u - A(u)\|\ls\|u - B(u)\|\ls 2\|u - A(u)\|.\]
On the other hand, we have for any \(u\in\widetilde{X}\),
\[\be{ll}
\langle J^{\prime}(u),u - B(u)\rangle&=\sum\limits_{i = 1}^{n(u)}\mu_{\lambda_{i}}(u)\langle J^{\prime}(u),u - A(z_{\lambda_{i}})\rangle\\
&=\sum\limits_{i = 1}^{n(u)}\mu_{\lambda_{i}}(u)\big(\langle J^{\prime}(u),u - A(u)\rangle+\langle J^{\prime}(u),A(u) - A(z_{\lambda_{i}})\rangle\big)\\
&\gl\sum\limits_{i = 1}^{n(u)}\mu_{\lambda_{i}}(u)\big(\langle T(u),u - A(u)\rangle - (\mathbf{f}(u),u - A(u))-\|J^{\prime}(u)\|\|A(u) - A(z_{\lambda_{i}})\|\big)\\
&\gl\sum\limits_{i = 1}^{n(u)}\mu_{\lambda_{i}}(u)\big(\langle T(u) - T(A(u)),u - A(u)\rangle - \|J^{\prime}(u)\|\|A(u) - A(z_{\lambda_{i}})\|\big)\\
&\gl\sum\limits_{i = 1}^{n(u)}\mu_{\lambda_{i}}(u)\big(c_{1}\|u - A(u)\|^{p} - \|J^{\prime}(u)\|\|A(u) - A(z_{\lambda_{i}})\|\big),
\en\eqno(2.7)\]
while by \((2.3)\)$\sim$\((2.5)\),\[\be{ll}\|J^{\prime}(u)\|\|A(u) - A(z_{\lambda_{i}})\|&\ls 2\|J^{\prime}(u_{\la_i})\|\cdot\frac{c_1}{8}(\|J^{\prime}(u_{\la_i})\|)^{-1}\|u_{\lambda_{i}} - A(u_{\lambda_{i}})\|^{p}\\
&\ls\frac{c_{1}}{2}\|u - A(u)\|^{p}.\\
\en\eqno(2.8)\]
It follows from (2.7) and (2.8) that
\[\langle J^{\prime}(u),u - B(u)\rangle\gl\frac{c_{1}}{2}\sum_{i = 1}^{n(u)}\mu_{\lambda_{i}}(u)\cdot\|u - A(u)\|^{p}=\frac{c_{1}}{2}\|u - A(u)\|^{p}.\]
The proof is complete.
\medskip

Next, we consider the following initial value problem in Banach space.
\[\left\{\be{ll}
\di\frac{du}{dt}=B(u(t)) - u(t)\\
u(0)=u_{0}.
\en\right.
\eqno(2.9)\]
Since \(u - B(u)\) is locally Lipschitz continuous, according to the theory of ordinary differential equations in Banach space, \((2.9)\) has a right saturated solution \(u(t,u_{0})\) for \(t\in[0,\eta(u_{0}))\), where \(0<\eta(u_{0})\ls+\infty\).
For any given \(t\in[0,\eta(u_{0}))\), we have (we simply denote \(u(t)=u(t,u_{0})\) below)
\[\be{ll}\frac{dJ(u(t))}{dt}&=\langle J^{\prime}(u(t)),u^{\prime}(t)\rangle=-\langle J^{\prime}(u(t)),u(t) -B(u(t))\rangle\\
&\ls-\frac{c_{1}}{2}\|u(t) - A(u(t))\|^{p}<0.
\en\eqno(2.10)\]
This shows that \(J(u(t))\) is decreasing in \(t\in[0,\eta(u_{0}))\). Hence, as in [10] we have the following definition.
\medskip

{\bf Definition 2.4}\
Let \(D\subset X\). If for any given \(u_{0}\in D\backslash K\), we have
\[\{u(t,u_{0}):t\in[0,\eta(u_{0}))\}\subset D,\]
then \(D\) is called an invariant set of  descending flow of \((2.9)\).
\medskip

{\bf Lemma 2.4}\
 Let \(D\subset X\) be a closed invariant set of  descending flow of \((2.9)\), \(d_{0}=\inf_{u\in D}J(u)>-\infty\). Then \(d_{0}\) is a critical value of \(J\), and there exists \(v_0\in D\cap K\) such that \(J(v_0)=d_{0}\).

{\bf Proof.}\ To show this lemma we will follows some ideas in [Theorem 2.1, 10].
By the definition of \(d_{0}\), for any given \(n\in \mb N\), there exists \(\overline{u}_{n}\in D\bk K\) such that \(J(\overline{u}_{n})\ls d_{0}+\frac{1}{n}\).
Consider the following initial value problem
\[\left\{\be{ll}
\di\frac{du}{dt}=B(u(t)) - u(t),\\
u(0)=\overline u_{n}.
\en\right.
\eqno(2.11)\]
Let the right saturated solution of \((2.11)\) be simply denoted as \(u_{n}(t)\) for $t\in[0,\eta(\overline{u}_{n}))$.  Similar to \((2.10)\), we have
\[\frac{dJ(u_{n}(t))}{dt}\ls-\frac{c_{1}}{2}\|u_{n}(t)-A(u_{n}(t))\|^{p},\ t\in[0,\eta(\overline{u}_{n})).\eqno(2.12)\]
Then we have the following two cases.

\((i)\) \(\eta(\overline{u}_{n})<\infty\). In this case, for \(\{t_{k}\}\subset[0,\eta(\overline{u}_{n}))\) being an increasing sequence, \(t_{k}\to\eta^-(\overline{u}_{n})\) (as \(k\to\infty\)), we have
\begin{align*}
\|u_{n}(t_{k})-u_{n}(t_{k-1})\|&\ls\int_{t_{k-1}}^{t_{k}}\|u_{n}(t)-B(u_{n}(t))\|dt\ls 2\int_{t_{k-1}}^{t_{k}}\|u_{n}(t)-A(u_{n}(t))\|dt\\
&\ls 2\left(\int_{t_{k-1}}^{t_{k}}\|u_{n}(t)-A(u_{n}(t))\|^{p}dt\right)^{\frac{1}{p}}(t_{k}-t_{k-1})^{\frac{1}{q}}\\
&\ls 2\left(\int_{t_{k-1}}^{t_{k}}\left(-\frac{2}{c_1}\frac{dJ(u_{n}(t))}{dt}\right)dt\right)^{\frac{1}{p}}(t_{k}-t_{k-1})^{\frac{1}{q}}\\
&=\frac{2^{\frac{1}{p}+1}}{c_1^{\frac{1}{p}}}\left(J(u_{n}(t_{k-1})-J(u_{n}(t_{k}))\right)^{\frac{1}{p}}(t_{k}-t_{k-1})^{\frac{1}{q}}\\
&\ls\frac{2^{\frac{1}{p}+1}}{c_1^{\frac{1}{p}}}(d_0+\frac{1}{n}-d_0)^{\frac{1}{p}}(t_{k}-t_{k-1})^{\frac{1}{q}}\\
&\ls\frac{2^{\frac{1}{p}+1}}{c_1^{\frac{1}{p}}}(t_{k}-t_{k-1})^{\frac{1}{q}},
\end{align*}
where \(\frac{1}{p}+\frac{1}{q}=1\), which shows that \(\{u_{n}(t_{k})\}_{k=1}^{\infty}\subset X\) is a Cauchy sequence, so there exists \(v_{n}\in X\) such that
\(u_{n}(t_{k})\to v_{n}\) in $X$ as \(k\to\infty\). According to the theory of ordinary differential equations in abstract space, \(v_{n}\in K\).

$(ii)$ $\eta(\overline{u}_{n}) = +\infty$. In this case, from \((2.12)\), we know that
\[
\frac{c_1}{2}\di\int_{0}^{+\infty}\|u_{n}(t)-A(u_{n}(t))\|^{p}dt\ls\di\int_{0}^{+\infty}\left(-\frac{dJ(u_{n}(t))}{dt}\right)dt
 = J(\overline{u}_{n})-J(u_{n}(t_{k}))<+\infty,
\]
and so \[\int_{0}^{+\infty}\|u_{n}(t)-A(u_{n}(t))\|^{p}dt<+\infty.\]
Consequently, there exists \(\{t_{k}\}_{k = 1}^{\infty}\subset[0,+\infty)\) such that \(t_{k}\to+\infty\), and \[\|u_{n}(t_{k})-A(u_{n}(t_{k}))\|\to 0\  \mbox{as}\  k\to\infty. \eqno(2.13)\]
Note that \(d_0\ls J(u_{n}(t_{k}))\ls d_0+\frac{1}{n}\) By (2.13) and (H$_4)$ we know that \(\{u_{n}(t_{k})\}_{k = 1}^{\infty}\) is bounded. Since \(X\) is a reflexive space, by  Alaoglu Theorem,
there exists a subsequence of \(\{u_{n}(t_{k})\}_{k = 1}^{\infty}\) (without loss of generality, we assume it is itself) and \(v_{n}\in X\) such that \(u_{n}(t_{k})\rightharpoonup v_{n}\) in $X$ as \(k\to\infty\).
Consequently,  we have
\[\langle J^{\prime}(v_{n}),u_{n}(t_{k})-v_{n}\rangle= o(1)\quad\text{as }k\to\infty,\eqno(2.14)\]
Since
\begin{align*}
\vert\langle J^{\prime}(u_{n}(t_{k}))-J^{\prime}(v_{n}),u_{n}(t_{k})-v_{n}\rangle\vert&\ls\vert\langle J^{\prime}(u_{n}(t_{k})),u_{n}(t_{k})-v_{n}\rangle\vert+\vert\langle J^{\prime}(v_{n}),u_{n}(t_{k})-v_{n}\rangle\vert\\
&\ls\vert J^{\prime}(u_{n}(t_{k}))\vert(\|u_{n}(t_{k})\|+\|v_{n}\|)+\vert\langle J^{\prime}(v_{n}),u_{n}(t_{k})-v_{n}\rangle\vert,
\end{align*}
by (2.14) we have $$\langle J^{\prime}(u_{n}(t_{k}))-J^{\prime}(v_{n}),u_{n}(t_{k})-v_{n}\rangle=o(1)\ \mbox{as}\ k\ri\infty.$$
On the other hand, we have
\[\be{ll}
\langle J^{\prime}(u_{n}(t_{k}))-J^{\prime}(v_{n}),u_{n}(t_{k})-v_{n}\rangle&=\langle T(u_{n}(t_{k}))-T(v_{n}),u_{n}(t_{k})-v_{n}\rangle\\
&+( \mathbf{f}(v_{n})-\mathbf{f}(u_{n}(t_{k})),u_{n}(t_{k})-v_{n}).
\en\eqno(2.15)\]
Noting that \(\mathbf{f}:E\to E^*\) is continuous and $X\hookrightarrow\hookrightarrow E$, so \(\mathbf{f}(u_{n}(t_{k}))\to \mathbf{f}(v_{n})\) as \(k\to\infty\) in \(E^*\). Then
\[
( \mathbf{f}(v_{n})-\mathbf{f}(u_{n}(t_{k})),u_{n}(t_{k})-v_{n})= o(1)\quad\text{as }k\to\infty.\eqno(2.16)
\]
From \((2.14)\sim(2.16)\), we know that as \(k\to\infty\)
\[
c_1\|u_{n}(t_{k})-v_{n}\|^{p}\ls\langle T(u_{n}(t_{k}))-T(v_{n}),u_{n}(t_{k})-v_{n}\rangle = o(1).
\]
Therefore, \(u_{n}(t_{k})\ri v_{n}\) in \(X\) as $k\ri\infty$, and so $v_n\in D$. By \((2.13)\) and (H$_4)$, we know that \(J^{\prime}(v_{n}) = 0\), and obviously
\(d_0\ls J(v_{n})\ls d_0+\frac{1}{n}\).

Both in the case (i) and (ii),  by the condition (H\(_{4})\), we know that \(\{v_{n}\}\) is bounded. Suppose \(v_{n}\rightharpoonup v_{0}\) in $X$. Similar to the above proof, we can show that \(v_{n}\to v_{0}\) in $X$. Then, we know that \(J(v_{0}) = d_0\), and obviously \(J^{\prime}(v_{0}) = 0\) and $v_0\in D$. The proof is complete.
\medskip

{\bf  Definition 2.5([10])}\  Let \(M\subset X\) be a connected invariant set of descending flow  of \((2.9)\), and \(D\subset M\) be a invariant set of  descending flow of \((2.9)\). Let
\[C_{M}(D)=\{v_{0}:v_{0}\in D \text{ or } v_{0}\in M\backslash D \text{ and there exists } t\in[0,\eta(v_{0})) \text{ such that }
u(t,v_{0})\in D\}.\]
Then  \(C_{M}(D)\) is called the descending flow invariant set expanded from \(D\) in \(M\). If \(C_{M}(D)=D\), then  \(D\) is called a complete descending flow invariant set of \((2.9)\) in \(M\).
\medskip

{\bf Lemma 2.5([10])}\ Let \(M\subset X\) be a connected invariant set of  descending flow of \((2.9)\), and \(D\) be an open invariant set of  descending flow of \((2.9)\) in \(M\), then

(1)\ \(C_{M}(D)\) is an open  invariant set of descending flow in \(M\);

(2)\ \(\partial_{M}C_{M}(D)\) (the boundary of \(C_{M}(D)\) in \(M\)) is a complete invariant set of  descending flow of \((2.9)\) in \(M\), and
\[\inf\limits_{u\in\partial_{M}C_{M}(D)}J(u)\gl\inf\limits_{u\in\partial_{M}(D)}J(u).\]
\medskip

{\bf Lemma 2.6([10])}\  Let \(U\) be a bounded connected open set in \(\mb R^{N}\) containing \((0,0)\), then there exists a connected component \(\Gamma\) in \(\partial U\)
such that for any ray \(l\) starting from the origin, \(l\cap\Gamma\neq\emptyset\)
\medskip

{\bf Proof of Theorem 2.1}\ The proof can be divided into the following three steps:

{\bf Step 1.}\ We first prove that  both \(\overline{D}_1\) and \(D_1\) are  invariant sets of descending flow of \((2.9)\). In fact, for any given \(u_{0}\in\overline{D}_1\) and $\lambda\gl 0$ small, \[u_{0}+\lambda(B(u_{0})-u_{0})=(1-\lambda)u_{0}+\lambda B(u_{0})\in\overline{D}_1.\]
Then, by the Brezis-Martin theorem (cf. Theorem 1.6.3 in [18]), we know that \(u(t,u_{0})\in\overline{D}_1\) for any given \(t\in[0,\eta(u_{0}))\). This implies that \(\overline{D}_1\) is a invariant set of  descending flow of \((2.9)\).
Furthermore, from \((2.9)\), we know that \(u(t,u_{0})\) satisfies the formula:
\[u(t,u_{0}) = e^{-t}u_{0}+\int_{0}^{t}e^{-(t - s)}B(u(s,u_{0}))ds,\ 0\ls t<\eta(u_{0}).\]
Therefore, for any given \(0\ls t<\eta(u_{0})\), we have
\[u(t,u_{0}) = e^{-t}u_{0}+(1 - e^{-t})\lim_{n\to\infty}\frac{1}{n}\sum_{k = 1}^{n}B\left(u(\ln(1 + \frac{k}{n}(e^{t}-1)),u_{0})\right).\eqno (2.17)\]
When \(u_{0}\in D_1\),
since for any given \(1\ls k\ls n\), \(u\left(\ln(1 + \frac{k}{n}(e^{t}-1)),u_{0}\right)\in\overline{D}_1\), thus
\[B\left(u(\ln(1 + \frac{k}{n}(e^{t}-1)),u_{0})\right)\in D^{+}(\frac{1}{2}\mu).\]
Noting that \(\overline{D}^{+}(\frac{1}{2}\mu)\) is a closed convex set, so
\[\lim_{n\to\infty}\frac{1}{n}\sum_{k = 1}^{n}B\left(u(\ln(1 + \frac{k}{n}(e^{t}-1)),u_{0})\right)\in\overline{D}^{+}(\frac{1}{2}\mu).\]
Consequently, by (2.17), for any given \(t\in[0,\eta(u_{0}))\),  \(u(t,u_{0})\in D_1\) for any \(u_{0}\in D_1\).
Therefore, \(D_1\) is a  invariant set of  descending flow of \((2.9)\).

Similarly, \(\overline{D}_2\) and \(D_2\) are  invariant sets of  descending flow of \((2.9)\).
 So, $D_3$ is an open invariant set of descending flow of (2.9).

{\bf Step 2.}  In this step we prove the existence of positive and negative solutions. We claim that $J$ is bounded from below on $\ov D_3$.  In fact, for each $u\in D_3$, we have $d(u, P)\ls \mu$ and $d(u,-P)\ls \mu$. In a reflexive Banach space, the distance from any point to a closed convex cone is attained. Thus we may take $w_1\in P$ and $w_2\in -P$ such that $d(u, P)=\|u-w_1\|$ and $d(u, -P)=\|u-w_2\|$. Write $\xi_u=u-w_1$ and $\eta_u=u-w_2$. Then we have $w_1-w_2=\eta_u-\xi_u$. Since $X\hookrightarrow\hookrightarrow E$, we have
 $$\|\xi_u\|_1\ls c_0\|\xi_u\|\ls c_0\mu,\  \|\eta_u\|_1\ls c_0\|\eta_u\|\ls c_0\mu.$$
Note that the norm of  $E$ is a Riesz norm and $0\ls w_1\ls w_1-w_2$. Then, we have
$$\|w_1\|_1\ls \|w_1-w_2\|_1=\|\eta_u-\xi_u\|_1\ls \|\eta_u\|_1+\|\xi_u\|_1\ls 2 c_0\mu,$$
and so $$\|u\|_1\ls \|w_1\|_1+\|\xi_u\|_1\ls 3c_0\mu.$$
Consequently, by (H$_2)$ we have for each $u\in D_3$,
$$|F(u)|=\left|\int_0^1( \mathbf{f}(su), u) ds\right|\ls \int_0^1g (s\|u\|_1)\|u\|_1 ds\ls 3c_0\mu g(3c_0\mu).$$
Thus,  we have for each $u\in D_3$,
$$J(u)=\int_0^1\langle T(su), u\rangle ds-\int_0^1( \mathbf{f}(su), u) ds\gl \frac{c_1}{p}\|u\|^p- 3c_0\mu g(3c_0\mu)\gl -3c_0\mu g(3c_0\mu).$$
Hence, $J$ is bounded from below on $\ov D_3$.

Take $e_1\in P\bk\{0\}$ and $e_2\in X\bk (D_1\cup D_2)$. Let $X_1=\mbox{span}\ \{e_1, e_2\}$. It follows from  (H$_3)$ that there exists  $R_0>0$ such that $J(u)<-3c_0\mu g(3c_0\mu)$ for each $u\in X_1$ with $\|u\|=R_0$.

 Clearly, \(D_3\cap\overline{D}_1\) is an open  invariant set of descending flow of (2.9) in \(\overline{D}_1\). It is easy to see that \(\partial_{\overline{D}_1}C_{\overline{D}_1}(D_3\cap\overline{D}_1)\neq\emptyset\). By Lemma 2.5,
\(\partial_{\overline{D}_1}C_{\overline{D}_1}(D_3\cap\overline{D}_1)\) is a complete descending flow invariant set of \((2.9)\) in \(\overline{D}_1\), and
$$c_{+}:=\inf\limits_{u\in\partial_{\overline{D}_1}C_{\overline{D}_1}(D_3\cap\overline{D}_1)} J(u)
\gl\inf\limits_{u\in\partial_{\overline{D}_1}(D_3\cap\overline{D}_1)}J(u)\gl-3c_0\mu g(3c_0\mu).
$$
Using Lemma 2.4, we know that \(c_{+}\) is a  critical value of \(J\), and there exists \[u_{1}\in\partial_{\overline{D}_1}C_{\overline{D}_1}(D_3\cap\overline{D}_1)\cap K\] such that \(J(u_{1}) = c_{+}\). Since $J$ has no critical point on $\big(\overline{D}_1\cup\overline{D}_1\big)\bk (P\cup (-P))$,  \(u_{1}\) is a positive critical point of \(J\), and thus, \(u_{1}\) is a positive solution of \((1.1)\). Similarly, \((1.1)\) has a negative solution \(u_{2}\) such that $J(u_2)=c_-$, where $$c_{-}:=\inf\limits_{u\in\partial_{\overline{D}_1}C_{\overline{D}_1}(D_3\cap\overline{D}_1)} J(u).$$

{\bf Step 3.}\ Finally, we prove that \((1.1)\) has at least one sign-changing solution. By Lemma 2.5, \(C_{X}(D_3)\) is an open set in \(X\).
Consequently, \(C_{X}(D_3)\cap X_{1}\) is an open set in \(X_{1}\).  It follows from the proof above that \(C_{X}(D_3)\cap X_{1}\subset B(0,R_{0})\). This indicates that
\(C_{X}(D_3)\cap X_{1}\) is a bounded open set in \(X_{1}\) containing the origin of $X_1$. By Lemma 2.6, there exists a connected component \(\Gamma_0\) of \(\partial(D_3\cap X_{1})\)
such that  \(l\cap\Gamma_0\neq\emptyset\) for any ray \(l\) starting from the origin of $X_1$. Let \({\Gamma}\) be the connected component of \(\partial_{X}C_{X}(D_3)\)
containing \(\Gamma_0\), then both \({\Gamma}\cap D_1\) and \({\Gamma}\cap D_2\) are  invariant sets of descending flow  of \((2.9)\) in $\Gamma$.
By Lemma 2.5, \(C_{{\Gamma}}({\Gamma}\cap D_1)\) and \(C_{{\Gamma}}({\Gamma}\cap D_2)\) are both non-empty open
 invariant sets of descending flow of  \((2.9)\). Since \({\Gamma}\) is connected, so
\[{\Gamma}_{1}:={\Gamma}\backslash(C_{{\Gamma}}({\Gamma}\cap D_1)\cup C_{{\Gamma}}({\Gamma}\cap D_2))\neq\emptyset.\]
Obviously, \({\Gamma}_{1}\) is an invariant set of  descending flow of \((2.9)\), and
\[\widetilde{c}:=\inf\limits_{u\in{\Gamma}_{1}}J(u)\gl\inf\limits_{u\in\partial D_3}J(u)\gl -3c_0\mu g(3c_0\mu).\]
Using Lemma 2.4, we know that \(\widetilde{c}\) is a critical value of \(J\), and there exists \(u_{3}\in{\Gamma}_{1}\) such that \(J(u_{3})=\widetilde{c}\).
It is easy to see that \(u_{3}\) is a sign-changing critical point of \(J\). Therefore, \((1.1)\) has a sign-changing solution \(u_{3}\). The proof is complete.
\medskip

{\bf Remark 2.2}\ If \(X\) is a Hilbert space, $(u, u^\pm)=(u^\pm, u^\pm)$ for any $u\in X$ and \(T = I\) is the identity mapping, then \(T\) satisfies \((\mbox{H}_{1})\) naturally.
\medskip

Under the assumptions of Theorem 2.1, we can further compare the critical values of sign-changing solutions, positive solutions and negative solutions.

{\bf Corollary 2.1}\ Suppose that the assumptions of Theorem 2.1 hold. Then there exists a  positive solution $x_1$, a negative solution $x_2$  and a sign-changing solution $x_3$, satisfying
\[
J(x_3)\gl \max\big\{J(x_1),J(x_2)\big\}.
\]

{\bf Proof.}\ Let $\Gamma$ be the set defined in  step 3 of  the proof of Theorem 2.1. From the proof of Theorem 2.1 above, $\ov {D}_1\cap \Gamma$ and $\ov {D}_2\cap \Gamma$  are closed invariant sets.  Then by Lemma 2.4, there exist critical points $x_1\in \ov{D}_1\cap \Gamma$ and $x_2\in \ov{D}_2\cap \Gamma$  such that
\[
J(x_1)=\inf_{u\in \ov{D}_1\cap \Gamma} J(u), J(x_2)=\inf_{u\in  \overline{D}_2\cap\Gamma} J(u).
\]
Since $x\in K$ implies $x=A(x)$, Lemma 2.3 yields $x_1\in P$, $x_2\in -P$.

 By Lemma 2.5,
\[
\overline{c}_3=\inf_{u\in \partial_{\Gamma}C_\Gamma({D}_1\cap \Gamma)} J(u)
\gl \inf_{u\in D_3} J(u)>-\infty
\]
is a critical value of $J$, and there exists a critical point $\overline{x}_3\in \partial_{\Gamma}C_\Gamma({D}_1\cap \Gamma)$ such that $J(\overline{x}_3)=\overline{c}_3$.
Take a sequence $\{y_n\}\subset C_{\Gamma}({D}_1\cap \Gamma)$ with $y_n\to \overline{x}_3$ as $n\to\infty$.

Consider the initial value problem
\[
\begin{cases}
\dfrac{du}{dt}=B(u(t))-u(t),\\[4pt]
u(0)=y_n,
\end{cases}
\]
which has a maximal global solution $u(t,y_n)$, with its existence interval  $[0,\eta(y_n))$. Then, there exists $t'>0$ such that $u(t',y_n)\in \ov{D}_1\cap \Gamma$, and so
\[
J(y_n)\gl J\big(u(t',y_n)\big)\gl J({x}_1).
\]
Letting $n\to\infty$, we obtain $\overline{c}_3\gl J({x}_1)$.

Similarly, there exists a critical point $\wi{x}_3\in \partial_{\Gamma}C_\Gamma({D}_2\cap \Gamma)$ such that
$$J(\wi{x}_3)=\wi{c}_3:=\inf_{u\in \partial_{\Gamma}C_\Gamma({D}_2\cap \Gamma)} J(u).$$
Here, we  may have $\ov x_3=\wi x_3$. Take $x_3$ be $\ov x_3$ or $\wi x_3$ such that $J(x_3)=\max\{\ov c_3,\wi c_3\}$. Then the conclusion holds. The proof is complete.

\section{Applications to Differential Equations Boundary Value Problems}
As an application of Theorem 2.1, we first consider the following \((p,2)\)-Laplacian equation boundary value problem
\[\left\{\be{ll}
-\Delta_{p}^{a_1}u(x)-\Delta_{2}^{a_2}u(x)=g(x,u)&\quad \text{in }\Omega,\\
u|_{\partial\Omega}=0,&\quad 2<p<p^*,
\en\right.\eqno(3.1)\]
where \(\Omega\subset \mb R^{N}\) is a bounded domain with a smooth boundary, and let \(C^{0,1}(\overline{\Omega})\) denote all Lipschitz
continuous functions. For \(u\in C^{0,1}(\overline{\Omega})\) and \(r\in(1,+\infty)\), \(\Delta_{r}^{u}\) we denote the weighted
\(r-\)Laplacian differential operator, defined as
\[\Delta_{r}^{u}(u)=\text{div}(a(x)|Du|^{r - 2}Du)\ \mbox{for all}\ u\in W_{0}^{1,r}(\Omega).\]
Here we assume that \(a_{1},a_{2}\in C^{1}(\overline{\Omega})\), and \(a_{1}(x),a_{2}(x)\gl\hat{c}>0\) for all \(x\in\overline{\Omega}\).
\medskip

 Concerning (3.1) we make the following assumptions:
\begin{itemize}
\item[($A_{1}$)]  There exists \(\hat{a}\in L^{\infty}(\Omega)\), \(p<\gamma<p^{*}\), such that
\(|g(x,u)|\ls\hat{a}(x)(1 + |u|^{\gamma-1})\), a.e. \(x\in\Omega\) and \(u\in\mathbb{R}\).
\item[($A_{2}$)] There are \(\gamma_{1}>p\) and \(M\gl 0\) such that, if \(G(x,u)=\int_{0}^{u}g(x,s)ds\), then \[
0<\gamma_1 G(x,s)\ls g(x,s)s\ \mbox{for a.a.}\ x\in \Om, \ \mbox{all}\ s\gl M,\]
\item[($A_{3}$)] There is $p<\tau<p^*$ such that \(\lim\limits_{u\to 0}\frac{g(x,u)}{|u|^{\tau-2}u}=0\) uniformly for a.a. \(x\in\Omega\).
\item[($A_{4}$)]  \(g(x,u)u\gl0\) a.a. \(x\in\Omega\) and \(u\in\mathbb{R}\).
\end{itemize}
\medskip
 Let $X=W_{0}^{1,p}(\Omega)$ and $E=L^{p^*}(\Om)$ with the norms
$$\|u\|=\big(\int_\Om |Du|^pdx\big)^{\frac{1}{p}}\ \mbox{for}\ u\in X,\ \|u\|_{p^*}=\big(\int_\Om |u|^{p^*}dx\big)^{\frac{1}{p^*}}\ \mbox{for}\ u\in E $$
respectively.  Let
$$P_1=\{u\in E: u(x)\gl 0 \ \mbox{a.e. in}\ \Om\}.$$
Then $P_1$ is cone in $E$ and $E$ is a Banach lattice with the ordering induced by the cone $P_1$. Moreover, $\|\cdot\|_1$ is a Riesz norm of $E$. Clearly, $X\hookrightarrow\hookrightarrow E$.

Let \(A_{p}^{a_{1}}:W_{0}^{1,p}(\Omega)\to W^{-1,p'}(\Omega)=W_{0}^{1,p'}(\Omega)^{*}\) and \(A_{2}^{a_{2}}:W_{0}^{1,2}(\Omega)\to W^{-1,2}(\Omega)=W_{0}^{1,2}(\Omega)^{*}\) (\(\frac{1}{p}+\frac{1}{p'}=1\)), be the nonlinear operators defined by
\[\langle A_{p}^{a_{1}}(u),h\rangle=\int_{\Omega}a_{1}(x)|Du|^{p - 2}(Du,Dh)_{\mathbb{R}^{N}}dx\ \mbox{for}\ u,h\in W_{0}^{1,p}(\Omega),\]
\[\langle A_{2}^{a_{2}}(u),h\rangle=\int_{\Omega}a_{1}(x)(Du,Dh)_{\mathbb{R}^{N}}dx\ \mbox{for}\ u,h\in W_{0}^{1,p}(\Omega).\]
Let
\(T(u)=A_{p}^{a_{1}}(u)+A_{2}^{a_{2}}(u)\) for $u\in W^{1,p}(\Omega)$. For each $u\in X$, let $\mathbf{f}: X\ri E^*$ be defined by $\mathbf{f}(u)(x)= g(x, u(x))$ a.a. $x\in \Om$. Then we have
\[(\mathbf{f}(u), w)=\int_\Om g(x, u(x)) w(x) dx\ \mbox{for any}\  w\in E.\]
We consider the \(C^{1}-\)functional \(J:W_{0}^{1,p}(\Omega)\to\mathbb{R}\) defined by
\[J(u)=\frac{1}{p}\int_{\Omega}a_{1}(x)|Du|^{p}dx+\frac{1}{2}\int_{\Omega}a_{2}(x)|Du|^{2}dx-\int_{\Omega}G(x,u)dx\] for all \(u\in W_{0}^{1,p}(\Omega)\). By a direct calculation, we know
\[\be{ll}\langle J^{\prime}(u),w\rangle&=\di\int_{\Omega}a_{1}(x)|Du|^{p - 2}(Du,Dw)_{\mathbb{R}^{N}}dx+\di\int_{\Omega}a_{2}(x)(Du,Dw)_{\mathbb{R}^{N}}dx-\int_{\Omega}g(x,u)wdx\\
&=\langle T(u),w\rangle-(\mathbf{f}(u),w).\en \]
\medskip

The following inequality is well known; See [8].

{\bf Lemma 3.1}\ For any \(\xi,\zeta\in\mathbb{R}^{N}\) it holds:
\[(|\xi|^{p - 2}\xi-|\zeta|^{p - 2}\zeta)\cdot(\xi-\zeta)\gl 2^{1 - p}|\xi-\zeta|^{p}\ \mbox{if}\ p\gl 2;\]
\[||\xi|^{p - 2}\xi-|\zeta|^{p - 2}\zeta|\ls(p - 1)(|\xi| + |\zeta|)^{p - 2}|\xi-\zeta|\ \mbox{if}\ p\gl2.\]
\medskip

In what follows we let $d_1, d_2,\cdots$ to denote  positive constants.

{\bf Lemma 3.2}\  Condition (H\(_{1})\) holds.

{\bf Proof}\ For any \(u,v\in W_{0}^{1,p}(\Omega)\), we have
\[\be{ll}
\langle T(u)-T(v),u-v\rangle
&=\di\int_{\Omega}a_{1}(x)(|Du|^{p - 2}Du-|Dv|^{p - 2}Dv)D(u-v)dx\\
&+\di\int_{\Omega}a_{2}(x)(Du-Dv)D(u-v)dx\\
&\gl2^{1 - p}\hat{c}\di\int_{\Omega}|Du - Dv|^{p}dx+\hat{c}\di\int_{\Omega}|Du - Dv|^{2}dx\\
&\gl2^{1 - p}\hat{c}\|u - v\|^{p}.
\en\]
Therefore, \(T:W_{0}^{1,p}(\Omega)\to W^{-1,p'}(\Omega)\) is strongly monotone.

 Let $\{u_k\}_{k=1}^\infty\subset W_{0}^{1,p}(\Omega)$ and $u\in W_{0}^{1,p}(\Omega)$ such that $u_k\ri u$ as \(k\to\infty\). Since \(W_{0}^{1,p}(\Omega)\hookrightarrow W_{0}^{1,2}(\Omega)\), we have
\[\be{ll}
\langle T(u_{k})-T(u),v\rangle&=\di\int_{\Omega}a_{1}(x)(|Du_{k}|^{p - 2}Du_{k}-|Du|^{p - 2}Du)\cdot Dvdx+\di\int_{\Omega}a_{2}(x)(Du_{k}-Du)\cdot Dvdx\\
&\ls\|a_{1}\|_{\infty}\Big(\di\int_{\Omega}| |Du_{k}|^{p - 2}Du_{k}-|Du|^{p - 2}Du|^{\frac{p}{p - 1}}dx\Big)^{\frac{p - 1}{p}}\Big(\di\int_{\Omega}|Dv|^{p}dx\Big)^{\frac{1}{p}}\\
&+\|a_{2}\|_{\infty}\Big(\di\int_{\Omega}| Du_{k}-Du|^{2}dx\Big)^{\frac{1}{2}}\Big(\di\int_{\Omega}|Dv|^{2}dx\Big)^{\frac{1}{2}}\\
&\ls\|a_{1}\|_{\infty}\Big((p-1)^{\frac{p}{p - 1}}\di\int_{\Omega}\Big(\big(|Du_{k}|+|Du|\big)^{p-2} |Du_{k}-Du|\Big)^{\frac{p}{p - 1}}dx\Big)^{\frac{p - 1}{p}}\big(\di\int_{\Omega}|Dv|^{p}dx\big)^{\frac{1}{p}}\\
&+\|a_{2}\|_{\infty}\Big(\di\int_{\Omega}| Du_{k}-Du|^{2}dx\Big)^{\frac{1}{2}}\Big(\di\int_{\Omega}|Dv|^{2}dx\Big)^{\frac{1}{2}}\\
&\ls d_1\Big(\di\int_{\Omega}\big(|Du_{k}|+|Du|)^{p}dx\Big)^{\frac{p-2}{p}}\Big(\di\int_{\Omega}|Du_{k}-Du|\big)^{p}dx\Big)^{\frac{1}{p}}\big(\di\int_{\Omega}|Dv|^{p}dx\big)^{\frac{1}{p}}\\
&+\|a_{2}\|_{\infty}\Big(\di\int_{\Omega}| Du_{k}-Du|^{2}dx\Big)^{\frac{1}{2}}\Big(\di\int_{\Omega}|Dv|^{2}dx\Big)^{\frac{1}{2}}.\\\en
\]
Consequently, we have as \(k\to\infty\),
\[\be{ll}
\|T(u_{k})-T(u)\|&=\sup\{\langle T(u_{k})-T(u),v\rangle:\|v\|=1\}\\
&\ls d_2(\|u_k\|+\|u\|)^{p-2}\|u_k-u\|+d_3\|u_k-u\|\to0.\en\eqno(3.2)\]
Therefore, \(T:W_{0}^{1,p}(\Omega)\to W^{-1,p^{\prime}}(\Omega)\) is continuous.

For each $u\in X$ we have \[\be{ll}
 \langle T(u), u^{\pm}\rangle&=\di\int_{\Omega}a_{1}(x)|Du|^{p-2}(Du,Du^\pm)_{\mathbb{R}^{N}}dx+\di\int_{\Omega}a_{2}(x)(Du,Du^\pm)_{\mathbb{R}^{N}}dx\\
 &=\di\int_{\Omega}a_{1}(x)|Du^\pm|^{p-2}(Du^\pm,Du^\pm)_{\mathbb{R}^{N}}dx+\di\int_{\Omega}a_{2}(x)(Du^\pm,Du^\pm)_{\mathbb{R}^{N}}dx\\
 &=\langle T(u^\pm), u^{\pm}\rangle.\en\]

Clearly,   $\mathbf{f}: E\ri E^*$ is continuous. Hence, (H\(_{1})\) holds. The proof is complete.
\medskip

{\bf Lemma 3.3}\  Condition (H\(_{2})\) holds.

{\bf Proof}\  It follows from $({\mbox{A}}_2)$ that \(g(x,0)=0\) a.a. \(x\in\Omega\). Then $\mathbf{f}$ is orthogonally additive. It follows from $(A_{1})$ and $(A_{3})$ that \[|g(x,u)|\ls d_4(|u|^{\tau-1}+|u|^{\ga-1})\ \mbox{a.a.}\ x\in\Omega\ \mbox{and}\ u\in\mathbb{R}.\]
Then by using H\"{o}lder inequality we have $$(\mathbf{f}(u),w)\ls d_5(\|u\|_{\tau}^{\tau-1}\|w\|_{\tau}+\|u\|_\ga^{\ga-1}\|w\|_\ga).$$
Since $L^\ga(\Om)\hookrightarrow E$ and $L^\tau(\Om)\hookrightarrow E$, we have
$$(\mathbf{f}(u),w)\ls d_6(\|u\|_{p^*}^{\tau-1}+\|u\|_{p^*}^{\ga-1})\|w\|_{p^*}.$$
 Let $g(s)=d_6 (s^{\tau-1}+ s^{\ga-1})$ for any $s\in \mb R^+$. Then the condition (H\(_{2})\) holds. The proof is complete.
\medskip

{\bf Lemma 3.4}\  Condition (H\(_{3})\) holds.

{\bf Proof}\  By \((A_{2})\), there exist \(\wi M>0\) and  $a>0$ such that
\[G(x,u)\gl a|u|^{\gamma_{1}} -\wi M, \ \forall u\in\mathbb{R}, \mbox{a.e.}\ x\in\Omega.\]
Then, for each $u\in X\bk\{0\}$ and $t\gl 0$, we have
\[\be{ll}J(tu)&=\frac{t^p}{p}\di\int_{\Omega}a_{1}(x)|Du|^{p}dx+\frac{t^2}{2}\di\int_{\Omega}a_{2}(x)|Du|^{2}dx-\int_{\Omega}G(x,u)dx\\
&\ls\frac{t^{p}}{p}\di\int_{\Omega}a_{1}(x)|Du|^{p}dx+\frac{t^{2}}{2}\di\int_{\Omega}a_{2}(x)|Du|^{2}dx-at^{\gamma_{1}}\int_{\Omega}|u|^{\gamma_{1}}dx+\wi M|\Omega|.\en\]
Noting that \(\gamma_{1}>p>2\), thus \(J(tu)\to-\infty\) as \(t\to+\infty\). The proof is complete.
\medskip

{\bf Lemma 3.5}\  Condition (H$_4)$ holds.

{\bf Proof}\ Let $[a,b]\subset \mb R$, \(\{u_{n}\} \subset J^{-1}([a, b])\) and $u_n-A(u_n)\ri 0$. Then, by (A$_2)$  we have
\[\be{ll}b+\|J'(u_n)\|\|u_n\|&\gl J(u_n)-\frac{1}{\ga_1}\langle J'(u_n), u_n \rangle\\
&\gl (\frac{1}{p}-\frac{1}{\ga_1})\|u_n\|^p+ (\frac{1}{2}-\frac{1}{\ga_1})|Du_n|_2^2+\di\int_\Om\big(\frac{1}{\ga_1}g(x,u_n(x))-G(x, u_n(x))\big)dx\\
&\gl (\frac{1}{p}-\frac{1}{\ga_1})\|u_n\|^p-d_7|\Om|.
\en\eqno(3.3)\]
 It follows from the proof of (3.2) that
$$\be{ll}\|J'(u_n)\|&=\|T(u_n)-T(A(u_n))\|\\
&\ls \big(d_2(\|u_n\|+\|A(u_n)\|)^{p-2}+d_3\big)\|u_n-A(u_n)\|\\
&\ls \big(d_2(2\|u_n\|+\|u_n-A(u_n)\|)^{p-2}+d_3\big)\|u_n-A(u_n)\|.\en\eqno(3.4)$$
Assume that $\|u_n-A(u_n)\|\ls 1$ for $n=1,2,\cdots$. By (3.3) and (3.4) we have
$$\big(\frac{1}{p}-\frac{1}{\ga_1}\big)\|u_n\|^p\ls \big(d_2(2\|u_n\|+1)^{p-2}+d_3\big)+d_7|\Om|+b.$$
Then we easily see that $\{u_n\}$ is bounded in $X$. Then, by (3.4) we have $J'(u_n)\ri 0$ in $X^*$ as $n\ri\infty$. The proof is complete.
\medskip

By Lemma 3.2$\sim $3.5, and Theorem 2.1 we can easily get the following result.

{\bf Theorem 3.1}\  Suppose \(({\mbox{A}}_{1})\sim({\mbox{A}}_{4})\) hold. Then \((3.1)\) admits at least one positive, one negative and one sign-changing solution.
\medskip

{\bf Remark 3.1}\ Recently, there are some papers studied the sign-changing solutions of $(p,q)$-Laplacian equations; see [16,17]. However, we need to point out that the methods in [16,17] are different from ours. For the existence of sign-changing solutions of $p$-Laplacian equations by using the method of descending flow invariant set one can refer to [15].

\end{document}